\documentclass[letterpaper, 10 pt, conference]{ieeeconf}  

\IEEEoverridecommandlockouts                              

\usepackage{amsmath} 
\usepackage{amssymb}  

\title{\LARGE \bf
Optimal Control Problem for Discrete-Time Systems with Colored Multiplicative Noise*
}

\author{Hongdan Li, Juanjuan Xu and
        Huanshui Zhang
\thanks{*This work was supported by the National Natural Science Foundation of China (under Grants 61633014,
61573220, 61573221).}
\thanks{The authors are with School of Control Science and Engineering, Shandong University, Jingshi Road 73,
        Jinan, 250061, P. R. China.
        {\tt\small  hszhang@sdu.edu.cn.}}%
}

\usepackage[letterpaper,top=72pt,bottom=54pt,left=54pt,right=54pt]{geometry}
\begin{document}

\maketitle
\thispagestyle{empty}
\pagestyle{empty}

\begin{abstract}

The optimal control problem for discrete-time systems with colored multiplicative noise is discussed in this paper. The problem will be more difficult to deal with than the case of white noise due to the correlation of the adjoining state. By solving the forward and backward stochastic difference equations (FBSDEs), the necessary and sufficient conditions for the solvability of the optimal control problems in both delay-free and one-step input delay case are given.

\end{abstract}

\section{INTRODUCTION}

The linear quadratic control was first introduced by Kalman \cite{1Kalman:60} in 1960, and attracted many other researchers to study it, see \cite{2Bismut:76}, \cite{3Zhang:98}, \cite{4Qi:17}, \cite{04Hou:17}, \cite{7Ju:18} and references therein.  It is generally known that uncertainty exists universally in practical application so that a renewed problem for systems with stochastic uncertainties has received much attention \cite{08Liang:18}, \cite{8Rami:00}, \cite{9Gershon:01}, \cite{09Gao:17}, \cite{13Xu:18} following the pioneering work by Wonham \cite{05Wonham:68}. For example, \cite{8Rami:00} discussed the indefinite LQ control problem for the discrete time system with state and control dependent noise and established the equivalence between the well-posedness and the attainability of the LQ problem. \cite{9Gershon:01} considered the  $H_{\infty}$  optimal control and filtering problem for linear discrete-time systems with stochastic uncertainties for the finite-horizon case.

It is worth noting that  the aforementioned works  almost considered the system that the coefficient of state and/or control variables involving  one multiplicative noise, such as the system in \cite{10Zhang:15} as follows:
\begin{eqnarray}
x_{k+1}=(A_{0}+A_{1}\omega_{k})x_{k}+(B_{0}+B_{1}\omega_{k})u_{k-d},\label{f1}
\end{eqnarray}
in which $\omega_{k}$ is a white noise with zero mean and variance $\sigma^{2}$. However, when noises $\omega_{k}$, $\omega_{k-1}$ simultaneously involved in the coefficient of control $u_{k-d}$, i.e.,
 \begin{eqnarray}
 x_{k+1}&=&(A_{0}+A_{1}\omega_{k})x_{k}
 +(B_{0}+B_{1}\omega_{k}\nonumber\\
 &&+B_{2}\omega_{k-1})u_{k-d},\label{f2}
\end{eqnarray}
$u_{k}$ is the input control with delay $d>0$ and $A, \bar{A}, B, \bar{B}, \tilde{B}$ are constant matrices with compatible dimensions, it can be seen that the system state is correlated at adjoining times. In order to distinguish the system (\ref{f2}),  we call this phenomenon as colored multiplicative noise systems.  Actually, this phenomenon exists in many fields such as in engineering field, see \cite{16Kay:81}, \cite{17Kay:72}, \cite{18Bryson:65} and references therein. But due to the correlation of the adjoining state, the LQ problem for the colored noise system with input delay will be more complex to solve.

Recently, some substantial progress for the optimal LQ control has been made by proposing the approach of solving the forward and backward differential/difference equations (FBDEs, for short), see \cite{10Zhang:15}, \cite{11Zhang:17} for details. Inspired by these works, we considered the linear quadratic optimal control problem for discrete time colored multiplicative noise system with one-step input delay or without delay. The contributions of this paper are as follows. Firstly, a necessary and sufficient condition for the optimal control problem to admit a unique solution is proposed in terms of the maximum principle. Secondly, the optimal controller and the optimal cost are explicitly presented via a coupled Riccati equation which is derived from the solution to the forward (the state equation) and backward (the costate equation) stochastic difference equations (FBSDEs, for short), which are more difficult to solve due to the correlation of the adjoining state compared with that in \cite{10Zhang:15}. Finally, the non-homogeneous relationship  between the costate and the state is established.

The rest of this paper is organized as follows.  The solvability of the optimal control problem for the delay-free case is given in Section 2. As to the one-step input delay case, the result is presented in Section 3.  Conclusions will be given in Section 4.

\emph{Notations:} $R^{n}$ stands for the usual $n$-dimensional Euclidean space;
$I$ denotes an identity matrix with appropriate dimension;  The superscript
$'$ represents the matrix transpose;  Real symmetric matrix $A > 0$ (or $\geq 0$) implies that A is strictly positive definite (or positive semi-definite). $\{\Omega,\mathcal{F},P,\{\mathcal{F}_{k}\}_{k\geq0} \}$ represents a
complete probability space, with natural filtration $\{\mathcal{F}_{k}\}_{k\geq0}$
generated by $\{\omega_{0}, \omega_{1}, \cdots, \omega_{k}\}$ augmented by all the $P$-null sets. $E[\cdot|\mathcal{F}_{k}]$ means the conditional expectation with respect to $\mathcal{F}_{k}$ and $\mathcal{F}_{-1}$ is understood as $\{\emptyset,\Omega\}$.

\section{Delay-free Case }
Consider the discrete-time stochastic system without delay
\begin{eqnarray}
 x_{k+1}=A_{k}x_{k}+B_{k}u_{k}, \label{f3}
\end{eqnarray}
where
\begin{eqnarray*}
  &&A_{k}=A_{0}+A_{1}\omega_{k}, \\ &&B_{k}=B_{0}+B_{1}\omega_{k}+B_{2}\omega_{k-1},
\end{eqnarray*}
$w_{k}$ is a scalar random white noise with zero mean and variance $\sigma^{2}$, and $A_{0}, A_{1}, B_{0}, B_{1}, B_{2}$ are constant matrices with compatible dimensions. And the following cost function:
\begin{eqnarray}
J_{N}&=&E\Big\{\sum_{k=0}^{N}(x'_{k}Qx_{k}+u'_{k}Ru_{k})\nonumber\\
&&+x'_{N+1}P_{N+1}x_{N+1}\Big\},\label{f4}
\end{eqnarray}
where $Q$, $R$ and $P_{N+1}$ are positive semi-definite matrices.

\textbf{ Problem 1:} \ \ Find a $\mathcal{F}_{k-1}$ measurable $u_{k}$ such that (\ref{f4}) is minimized subject to (\ref{f3}).

By Pontryagin's maximum principle, it yields the following costate equations
\begin{eqnarray}
\lambda_{k-1}=E[A'_{k}\lambda_{k}|\mathcal{F}_{k-1}]+Qx_{k},\label{f5}
\end{eqnarray}
with the terminal value
\begin{eqnarray}
\lambda_{N}=P_{N+1}x_{N+1},\label{f6}
\end{eqnarray}
and the equilibrium condition
\begin{eqnarray}
0=Ru_{k}+E[B'_{k}\lambda_{k}|\mathcal{F}_{k-1}].\label{f7}
\end{eqnarray}

To facilitate the explanation of  Problem 1, we will introduce the following difference equation as
\begin{eqnarray}
P_{k}=A'_{0}P_{k+1}A_{0}+\sigma^{2}A'_{1}P_{k+1}A_{1}+Q
-M'_{k}\Upsilon_{k}^{-1}M_{k},\label{f8}
\end{eqnarray}
where
\begin{eqnarray}
\Upsilon_{k}
&=&R+\sigma^{2}B'_{1}P_{k+1}B_{1}+(B_{0}+B_{2}\omega_{k-1})'P_{k+1}\nonumber\\
&&\times (B_{0}+B_{2}\omega_{k-1}), \label{f9}\\
M_{k}&=&B'_{0}P_{k+1}A_{0}+\sigma^{2}B'_{1}P_{k+1}A_{1}\nonumber\\
&&+\omega'_{k-1}B'_{2}P_{k+1}A_{0}.\label{f10}
\end{eqnarray}

\textbf{Remark 1: } \ \ When $B_{2}=0$, the above equation can be reexpressed as
\begin{eqnarray}
P_{k}&=&A'_{0}P_{k+1}A_{0}+\sigma^{2}A'_{1}P_{k+1}A_{1}+Q\nonumber\\
&&-\bar{M}_{k}'R^{-1}_{k}\bar{M}_{k},\label{f11}
\end{eqnarray}
with
\begin{eqnarray}
\bar{M}_{k}&=&B'_{0}P_{k+1}A_{0}+\sigma^{2}B'_{1}P_{k+1}A_{1},\label{f12}\\
R_{k}&=&R+B'_{0}P_{k+1}B_{0}+\sigma^{2}B'_{1}P_{k+1}B_{1},\label{f13}
\end{eqnarray}
which is the Riccati equation of white noise case.

The following is the introduction of the main theorem in this section.

\textbf{Theorem 1:} \ \ Problem 1 has a unique solution if and only
if $\Upsilon_{k}>0$. In this case, the optimal controller $u_{k}$ is stated as
\begin{eqnarray}
u_{k}=-\Upsilon_{k}^{-1}M_{k}x_{k}.\label{f14}
\end{eqnarray}

The associated optimal value of (\ref{f3}) is given by
\begin{eqnarray}
J_{N}^{\ast}
&=&E(x_{0}'P_{0}x_{0}).\label{f15}
\end{eqnarray}
Moreover, the optimal costate $\lambda_{k-1}$ and state $x_{k}$ satisfy the following non-homogeneous
relationship
\begin{eqnarray}
\lambda_{k-1}=P_{k}x_{k}.\label{f16}
\end{eqnarray}

\textbf{Proof:}\ \ \emph{``Necessity"}:  Assume that Problem 1 admits a unique solution, we will adopt  induction to illustrate $\Upsilon_{k}>0$ in (\ref{f9}). For convenience, we show
\begin{eqnarray}
\bar{J}_{k}&\triangleq& E\Big[\sum_{i=k}^{N}(x_{i}'Qx_{i}+u_{i}'Ru_{i})\nonumber\\
&&+x_{N+1}'P_{N+1}x_{N+1}|\mathcal{F}_{k-1}\Big]. \label{f100}
\end{eqnarray}
 Firstly, for $k=N$, from (\ref{f3}), we know that $\bar{J}_{N}$ can be expressed as a quadratic function of $x_{N}$ and $u_{N}$.  Letting $x_{N}=0$,  then $\bar{J}_{N}=u_{N}'\Upsilon_{N}u_{N}$, for the uniqueness of solution to Problem 1, it is clear that the optimal controller is $u_{N}=0$ and the optimal cost is 0.
 Therefore, for nonzero controller $u_{N}$, we have
\begin{eqnarray}
\bar{J}_{N}=u_{N}'\Upsilon_{N}u_{N}>0, \label{f101}
\end{eqnarray}
so does $\Upsilon_{N}>0$.

In this case, from (\ref{f3}), (\ref{f6}) and (\ref{f7}), we have
\begin{eqnarray}
0&=&Ru_{N}+E[B'_{N}\lambda_{N}|\mathcal{F}_{N-1}]\nonumber\\
&=&Ru_{N}+E[(B_{0}+B_{1}\omega_{N}+B_{2}\omega_{N-1})'P_{N+1}(A_{0}\nonumber\\
&&+A_{1}\omega_{N})
x_{N}+(B_{0}+B_{1}\omega_{N}+B_{2}\omega_{N-1})'P_{N+1}\nonumber\\
&&\times(B_{0}+B_{1}\omega_{N}+B_{2}\omega_{N-1})u_{N}|\mathcal{F}_{N-1}]\nonumber\\
&=&(R+B'_{0}P_{N+1}B_{0}+B'_{0}P_{N+1}\bar{B}_{2}\omega_{N-1}
+\sigma^{2}B'_{1}\nonumber\\
&&\times P_{N+1}B_{1}+\omega'_{N-1}B'_{2}P_{N+1}B_{0}
+\omega'_{N-1}B'_{2}\nonumber\\
&&\times P_{N+1}B_{2}\omega_{N-1})u_{N}+(B'_{0}P_{N+1}A_{0}\nonumber\\
&&+\sigma^{2}A'_{1}P_{N+1}A_{1}
+\omega'_{N-1}B'_{2}P_{N+1}A_{0})x_{N}\nonumber\\
&=&\Upsilon_{N}u_{N}+M_{N}x_{N}.\label{f102}
\end{eqnarray}
Therefore, the optimal controller can be computed as
\begin{eqnarray}
u_{N}=-\Upsilon_{N}^{-1}M_{N}x_{N},\label{f17}
\end{eqnarray}
which is correspond to (\ref{f14}) with $k=N$.

As to $\lambda_{N-1}$, from (\ref{f3}), (\ref{f6}) and (\ref{f17}), we can obtain that
\begin{eqnarray}
\lambda_{N-1}&=&E[A'_{N}\lambda_{N}|\mathcal{F}_{N-1}]+Qx_{N}\nonumber\\
&=&E[(A_{0}+\bar{A}_{1}\omega_{N})'P_{N+1}(A_{0}+\bar{A}_{1}\omega_{N})x_{N}\nonumber\\
&&+(A_{0}+\bar{A}_{1}\omega_{N})'P_{N+1}(B_{0}+\bar{B}_{1}\omega_{N}
+\bar{B}_{2}\omega_{N-1})\nonumber\\
&&\times u_{N}|\mathcal{F}_{N-1}]+Qx_{N}\nonumber\\
&=&(A'_{0}P_{N+1}A_{0}+\sigma^{2}\bar{A}'_{1}P_{N+1}\bar{A}_{1}+Q)x_{N}\nonumber\\
&&+(A'_{0}P_{N+1}B_{0}+A'_{0}P_{N+1}\bar{B}_{2}\omega_{N-1}
\nonumber\\
&& +\sigma^{2}\bar{A}'_{1}P_{N+1}\bar{B}_{1})u_{N}\nonumber\\
&=&(A'_{0}P_{N+1}A_{0}+\sigma^{2}\bar{A}'_{1}P_{N+1}\bar{A}_{1}+Q\nonumber\\
&&-M'_{N}\Upsilon_{N}^{-1}M_{N})x_{N},\label{f18}
\end{eqnarray}
from (\ref{f6}), it is clear to see $\lambda_{N-1}$ has the same form with (\ref{f12}) with $k=N$.

In order to complete the proof by induction, we take any $m$ with $0\leq m\leq N$. And for all $k\geq m+1$, we make the following assumptions: first, in (\ref{f9}) $\Upsilon_{k}>0$; second, $\lambda_{k-1}=P_{k}x_{k}$ in which $P_{k}$ satisfies (\ref{f6})-(\ref{f8}); third, $u_{k}=-\Upsilon_{k}^{-1}M_{k}x_{k}$ is the optimal control.
In consideration of these assumptions, next we'll verify that these are all satisfied for $k=m$.

First of all, we will illustrate $\Upsilon_{m}>0$.
From (\ref{f3}), (\ref{f5}) and (\ref{f7}), it yields that
\begin{eqnarray}
E[x_{k}'\lambda_{k-1}-x_{k+1}'\lambda_{k}]
&=&E\Big\{x_{k}'E[A_{k}\lambda_{k}|\mathcal{F}_{k-1}]
+x_{k}'Qx_{k}\nonumber\\
&&-x_{k}'A_{k}\lambda_{k}-u_{k}'B_{k}'\lambda_{k}\Big\}\nonumber\\
&=&E[x_{k}'Qx_{k}-u_{k}'E(B_{k}'\lambda_{k}|\mathcal{F}_{k-1})]\nonumber\\
&=&E[x_{k}'Qx_{k}+u_{k}'Ru_{k}].\label{f103}
\end{eqnarray}
Now we put the above-mentioned formula count up from $k=m+1$ to $k=N$ on both sides, then
\begin{eqnarray}
&&E[x_{m+1}'\lambda_{m}-x_{N+1}'P_{N+1}x_{N+1}]\nonumber\\
&=&\sum_{k=m+1}^{N}E[x_{k}'Qx_{k}+u_{k}'Ru_{k}].\label{f100}
\end{eqnarray}
Therefore, we have
\begin{eqnarray}
\bar{J}_{m}&=& E\Big[\sum_{i=m+1}^{N}(x_{i}'Qx_{i}+u_{i}'Ru_{i})+x_{N+1}'P_{N+1}x_{N+1}
\nonumber\\
&&+(x_{m}'Qx_{m}+u_{m}'Ru_{m})|\mathcal{F}_{m-1}\Big]\nonumber\\
&=&E[x_{m+1}'\lambda_{m}+(x_{m}'Qx_{m}+u_{m}'Ru_{m})|\mathcal{F}_{m-1}]. \label{f105}
\end{eqnarray}
To check $\Upsilon_{m}$, let $x_{m}=0$, thus,
\begin{eqnarray}
\bar{J}_{m}&=&E[u_{m}'B_{m}'\lambda_{m}+u_{m}'Ru_{m}|\mathcal{F}_{m-1}]\nonumber\\
&=&E[u_{m}'B_{m}'P_{m+1}B_{m}u_{m}+u_{m}'Ru_{m}|\mathcal{F}_{m-1}]\nonumber\\
&=&u_{m}'\Upsilon_{m}u_{m}.\label{f106}
\end{eqnarray}
For any nonzero control $u_{m}$, in view of the uniqueness of optimal control, we can obtain that $\Upsilon_{m}>0$. Next we will show the expression of the optimal $u_{m}$. From (\ref{f3}), (\ref{f5}) and (\ref{f7}), we have
\begin{eqnarray}
0&=&Ru_{m}+E[B'_{m}\lambda_{m}|\mathcal{F}_{m-1}]\nonumber\\
&=&Ru_{m}+E[B'_{m}P_{m+1}x_{m+1}|\mathcal{F}_{m-1}]\nonumber\\
&=&Ru_{m}+E[(B_{0}+B_{1}\omega_{m}+B_{2}\omega_{m-1})'P_{m+1}\nonumber\\
&&\times (A_{0}+A_{1}\omega_{m})
x_{m}+(B_{0}+B_{1}\omega_{m}+B_{2}\omega_{m-1})'\nonumber\\
&&\times
P_{m+1}(B_{0}+B_{1}\omega_{m}+B_{2}\omega_{m-1})u_{m}|\mathcal{F}_{m-1}]\nonumber\\
&=&(R+B'_{0}P_{m+1}B_{0}+B'_{0}P_{m+1}B_{2}\omega_{m-1}\nonumber\\
&&
+\sigma^{2}B'_{1}P_{m+1}B_{1}+\omega'_{m-1}B'_{2}P_{m+1}B_{0}\nonumber\\
&&+\omega'_{m-1}B'_{2}P_{m+1}B_{2}\omega_{m-1})u_{m}+(B'_{0}P_{m+1}A_{0}\nonumber\\
&&+\sigma^{2}B'_{1}P_{m+1}A_{1}+\omega'_{m-1}B'_{2}P_{m+1}A_{0})x_{m},\label{f19}
\end{eqnarray}
thus, the optimal control $u_{m}$ can be obtained as
\begin{eqnarray}
u_{m}=-\Upsilon_{m}^{-1}M_{m}x_{m}, \label{f20}
\end{eqnarray}
in which $\Upsilon_{m}, M_{m}$ are defined as in (\ref{f9}) and (\ref{f10}).

Next we will investigate the relationship between costate and state in the case of $k=m$. Considering (\ref{f3}), (\ref{f5}) and (\ref{f20}), it yields that
 \begin{eqnarray}
\lambda_{m-1}&=&E[A'_{m}\lambda_{m}|\mathcal{F}_{m-1}]+Qx_{m}\nonumber\\
&=&E[(A_{0}+A_{1}\omega_{m})'P_{m+1}(A_{0}+A_{1}\omega_{m})x_{m}\nonumber\\
&&+(A_{0}+A_{1}\omega_{m})'P_{m+1}(B_{0}+B_{1}\omega_{m}\nonumber\\
&&+B_{2}\omega_{m-1})u_{m}|\mathcal{F}_{m-1}]+Qx_{m}\nonumber\\
&=&(A'_{0}P_{m+1}A_{0}+\sigma^{2}A'_{1}P_{m+1}A_{1}+Q)x_{m}\nonumber\\
&&+(A'_{0}P_{m+1}B_{0}+A'_{0}P_{m+1}B_{2}\omega_{m-1}\nonumber\\
&&+\sigma^{2}A'_{1}P_{m+1}B_{1})u_{m}\nonumber\\
&=&(A'_{0}P_{m+1}A_{0}+\sigma^{2}A'_{1}P_{m+1}A_{1}+Q\nonumber\\
&&-M'_{m}\Upsilon_{m}^{-1}M_{m})x_{m}.\label{f21}
\end{eqnarray}
Seeing that (\ref{f8})-(\ref{f10}), it implies that the aforementioned equation satisfies (\ref{f16}). By induction, we complete the Necessity proof.

\emph{``Sufficiency"}: When $\Upsilon_{k}>0$ is satisfied, the unique optimal controller of Problem 1 and the optimal cost functional will be illustrated, respectively.

In order to illustrate the main result more clearly, we first define a function $V(k)$ as
\begin{eqnarray}
V(k,x_{k})=E[x_{k}'P_{k}x_{k}],\label{f22}
\end{eqnarray}
in which $P_{k}$ is as in (\ref{f8}).

In consideration of (\ref{f3}) and (\ref{f8})-(\ref{f10}), we can easily calculate that
\begin{eqnarray}
&&V(k)-V(k+1)\nonumber\\
&=&E[x_{k}'P_{k}x_{k}-x_{k+1}'P_{k+1}x_{k+1}]\nonumber\\
&=&E[x_{k}'P_{k}x_{k}-x_{k}'A_{k}'P_{k+1}A_{k}x_{k}-x_{k}'A_{k}'P_{k+1}B_{k}u_{k}
\nonumber\\
&&-u_{k}'B_{k}'P_{k+1}A_{k}x_{k}-u_{k}'B_{k}'P_{k+1}B_{k}u_{k}]\nonumber\\
&=&E[x_{k}'P_{k}x_{k}-x_{k}'(A_{0}'P_{k+1}A_{0}+\sigma^{2}A_{1}'P_{k+1}A_{1})x_{k}\nonumber\\
&&-x_{k}'(A_{0}'P_{k+1}B_{0}+\sigma^{2}A_{1}'P_{k+1}B_{1}+A_{0}'P_{k+1}\nonumber\\
&&\times B_{2}\omega_{k-1})u_{k}-u_{k}'(B_{0}'P_{k+1}A_{0}+\sigma^{2}B_{1}'P_{k+1}A_{1}\nonumber\\
&&+\omega_{k-1}'B_{2}'P_{k+1}A_{0})x_{k}-u_{k}'(B_{0}'P_{k+1}B_{0}\nonumber\\
&&+\sigma^{2}B_{1}'P_{k+1}B_{1}+B_{0}'P_{k+1}B_{2}\omega_{k-1}\nonumber\\
&&+\omega_{k-1}'B_{2}'P_{k+1}B_{0}+\omega_{k-1}'B_{2}'P_{k+1}B_{2}\omega_{k-1})u_{k}]\nonumber\\
&=&E[x_{k}'(Q-M_{k}'\Upsilon_{k}^{-1}M_{k})x_{k}-u_{k}'(\Upsilon_{k}-R)u_{k}\nonumber\\
&&-u_{k}'M_{k}x_{k}-x_{k}'M_{k}'u_{k}]\nonumber\\
&=&E[x_{k}'Qx_{k}+u_{k}'Ru_{k}-(u_{k}+\Upsilon_{k}^{-1}M_{k}x_{k})'
\Upsilon_{k}(u_{k}\nonumber\\
&&+\Upsilon_{k}^{-1}M_{k}x_{k})].\label{f23}
\end{eqnarray}
Adding from $k=0$ to $k=N$ on both sides of the aforementioned equation, the cost functional (1.2) can be rewritten as
\begin{eqnarray}
J_{N}&=&E\Bigg[\sum_{i=0}^{N}[(u_{k}+\Upsilon_{k}^{-1}M_{k}x_{k})'\Upsilon_{k}
(u_{k}+\Upsilon_{k}^{-1}M_{k}x_{k})]\nonumber\\
&&+x_{0}'P_{0}x_{0}\Bigg] \nonumber\\
&\geq&E(x_{0}'P_{0}x_{0}).\label{f24}
\end{eqnarray}
The condition for the above inequalities to be established is  $\Upsilon_{k}>0$. Hence, it's easy to see that the optimal cost is $J_{N}^{\ast}=E(x_{0}'P_{0}x_{0})$. And the optimal control has the same form as in (\ref{f14}).  This sufficient proof is completed.

\section{ One-step Delay Case}
Considering the following system with one-step input delay:
\begin{eqnarray}
 x_{k+1}&=&(A_{0}+A_{1}\omega_{k})x_{k}
 +(B_{0}+B_{1}\omega_{k}\nonumber\\
 &&+B_{2}\omega_{k-1})u_{k-1}.\label{f25}
\end{eqnarray}
The cost functional is
 \begin{eqnarray}
J_{N}&=&E\Big\{\sum_{k=0}^{N}x'_{k}Qx_{k}+\sum_{k=1}^{N}u'_{k-1}Ru_{k-1}\nonumber\\
&&+x'_{N+1}P_{N+1}x_{N+1}\Big\},\label{f26}
\end{eqnarray}
where $Q$, $R$ and $P_{N+1}$ are positive semi-definite matrices.

\textbf{ Problem 2:} \ \ Find a $\mathcal{F}_{k-1}$ measurable $u_{k}$ such that (\ref{f26}) is minimized subject to (\ref{f24}).

By stochastic maximum principle, we can obtain the forward and backward stochastic difference equations with one-step input delay.
\begin{eqnarray}
\left\{
\begin{array}{llll}
\lambda_{k-1}=E[(A_{0}+A_{1}\omega_{k})'\lambda_{k}|\mathcal{F}_{k-1}]+Qx_{k},\\
\lambda_{N}=P_{N+1}x_{N+1},\\
0=Ru_{k-1}+E[(B_{0}+B_{1}\omega_{k}\\
+B_{2}\omega_{k-1})'\lambda_{k}|\mathcal{F}_{k-2}],\\
x_{k+1}=(A_{0}+A_{1}\omega_{k})x_{k}
 +(B_{0}+B_{1}\omega_{k}\\
 +B_{2}\omega_{k-1})u_{k-1},
\end{array}
\right.\label{f27}
\end{eqnarray}
with initial values $x_{0}$ and $u_{-1}$.

Similar to the derivation of the solution to Problem 1, we have the solution to Problem 2.

\textbf{Theorem 2:} \ \ Problem 2 is uniquely solvable if and only
if $R_{k}>0$. In this case, the optimal controller $u_{k}$ is stated as
\begin{eqnarray}
u_{k}=-R_{k+1}^{-1}[T^{0}_{k+1}x_{k}+(T^{1}_{k+1}
+F_{k+1}B_{2}\omega_{k-1})u_{k-1}],\label{f28}
\end{eqnarray}
in which
\begin{eqnarray}
\hspace{-3mm}R_{k}
\hspace{-3mm}&=&\hspace{-3mm}R+B'_{0}P_{k+1}B_{0}+\sigma B'_{1}P_{k+1}B_{1}+\sigma B'_{2}(P_{k+1}\nonumber\\
\hspace{-3mm}&&\hspace{-3mm}-F'_{k+1}R^{-1}_{k+1}F_{k+1})B_{2}-(T^{1}_{k+1})'R^{-1}_{k+1}T^{1}_{k+1}, \label{f29}\\
T^{0}_{k}\hspace{-3mm}&=&\hspace{-3mm}B'_{0}P_{k+1}A^{2}_{0}+\sigma B'_{1}P_{k+1}A_{1}A_{0}+\sigma B'_{2}(P_{k+1}A_{0}\nonumber\\
\hspace{-3mm}&&\hspace{-3mm}-F'_{k+1}R^{-1}_{k+1}T^{0}_{k+1})A_{1}
-(T^{1}_{k+1})'R^{-1}_{k+1}T^{0}_{k+1}A_{0},\label{f30}\\
T^{1}_{k}\hspace{-3mm}&=&\hspace{-3mm}B'_{0}P_{k+1}A_{0}B_{0}+\sigma B'_{1}P_{k+1}A_{1}B_{0}+\sigma B'_{2}(P_{k+1}A_{0}\nonumber\\
\hspace{-3mm}&&\hspace{-3mm}-F'_{k+1}R^{-1}_{k+1}T^{0}_{k+1})B_{1}
-(T^{1}_{k+1})'R^{-1}_{k+1}T^{0}_{k+1}B_{0},\label{f31}\\
F_{k}\hspace{-3mm}&=&\hspace{-3mm}B'_{0}P_{k+1}A_{0}+\sigma B'_{1}P_{k+1}A_{1}-(T^{1}_{k+1})'R^{-1}_{k+1}T^{0}_{k+1},\label{f32}
\end{eqnarray}
with
\begin{eqnarray}
P_{k}\hspace{-3mm}&=&\hspace{-3mm}A'_{0}P_{k+1}A_{0}+\sigma A'_{1}P_{k+1}A_{1}+Q\nonumber\\
\hspace{-3mm}&&\hspace{-3mm}-(T^{0}_{k+1})'R^{-1}_{k+1}T^{0}_{k+1},\label{f33}
\end{eqnarray}
with terminal value $P_{N+1}$.
Moreover, the optimal costate $\lambda_{k-1}$ and state $x_{k}$ satisfy the following non-homogeneous
relationship
\begin{eqnarray}
\lambda_{k-1}&=&P_{k}x_{k}-(F_{k}+M_{k}\omega_{k-1})'R^{-1}_{k}
[T^{0}_{k}x_{k-1}\nonumber\\
&&+(T^{1}_{k}
+F_{k}B_{2}\omega_{k-2})u_{k-2}].\label{f34}
\end{eqnarray}
\textbf{Proof:} Following the proof of Theorem 1, the above result can be similarly obtained, so we omit it here.

\section{CONCLUSIONS}

 This paper mainly studied the linear quadratic regulation  problem for discrete-time systems with colored multiplicative noise for both delay-free and one-step input delay case. The necessary and sufficient condition for the solvability of optimal control problem was presented by solving the FBSDEs derived from the maximum principle. Moreover, the optimal controller and cost were given in terms of the coupled difference equations developed in this paper. For any input delay $d>0$, the systems are more general but more difficult to deal with. Therefore, the optimal control problem for linear systems with any input delay and colored multiplicative noise are worth considering in the future.

\addtolength{\textheight}{-12cm}   





\end{document}